\documentclass[centertags,12pt]{amsart}
\usepackage{latexsym}
\usepackage{amssymb}
\numberwithin{equation}{section}
\makeatletter
\renewcommand{\subsection}{\@startsection
{subsection}{2}{0mm}{\baselineskip}{-0.25cm}
{\normalfont\normalsize\bf}}
\makeatother
\newtheorem{theorem}{Theorem}[section]
\newtheorem{proposition}[theorem]{Proposition}
\newtheorem{corollary}[theorem]{Corollary}
\newtheorem{lemma}[theorem]{Lemma}
\newtheorem*{fund}{Fundamental Linear Equivalence}
{\theoremstyle{remark} 
\newtheorem*{claim*}{Claim}}
\theoremstyle{definition}
\newtheorem{remark}[theorem]{Remark}
\newtheorem{example}[theorem]{Example}
\newtheorem{question}[theorem]{Question}
\def\F{\mathbb F}
\def\P{\mathbb P}

\def\Z{\mathbb Z}
\def\cD{\mathcal D}
\def\cL{\mathcal L}
\def\cK{\mathcal K}
\def\cW{\mathcal W}
\def\fq{\mathbb F_{q^2}}
\def\xq{X\! \mid\! \mathbb F_{q^2}}
\def\supp{{\rm Supp}}
\def\df{{\rm div}}
\def\dim{{\rm dim}}
\def\deg{{\rm deg}}
\def\frx{{\rm Fr}_{X}}
\begin{document}
\author[A.~Garcia]{Arnaldo Garcia}
\author[F.~Torres]{Fernando Torres}\thanks{1991 Math. Subj. Class.:
Primary 11G, Secondary 14G}\thanks{ Both authors were partially supported by
Cnpq-Brazil and by PRONEX \# 41.96.0883.00}
\title[Classical maximal curves]{On maximal curves having classical\\
Weierstrass gaps}
\address{IMPA, Est. Dona Castorina 110, Rio de Janeiro, 22460-320-RJ,   
Brazil}
\email{garcia@impa.br}
\address{IMECC-UNICAMP, Cx. P. 6065, Campinas, 13083-970-SP, Brazil}
\email{ftorres@ime.unicamp.br}
\begin{abstract}
We study geometrical properties of maximal curves having 
classical Weierstrass gaps. 
\end{abstract}
\maketitle
\section{Introduction}\label{1}
Let $X$ be a projective geometrically irreducible nonsingular algebraic
curve of genus $g$, defined over a finite field $\mathbb F_\ell$ with
$\ell$ elements. In 1948, A. Weil \cite{weil} proved the Riemann
hypothesis for curves over finite fields which states that
$$
|\# X(\mathbb F_\ell)-(\ell+1)|\le 2g\sqrt{\ell}\, ,
$$
where $X(\mathbb F_\ell)$ denotes the set of $\mathbb F_\ell$-rational
points of the curve $X$. This bound was proved for elliptic curves by
Hasse.

The curve $X$ is called {\it maximal over $\mathbb F_\ell$} (in this case
$\ell$ must be a square; say $\ell=q^2$) if $\# X(\fq)$ attains the 
Hasse-Weil upper bound, that is
$$
\#X(\fq)=q^2+1+2gq\, .
$$
The genus of a maximal curve over $\fq$ satisfies \cite{ihara},
\cite{stix}, \cite{ft}
$$
g\le \frac{(q-1)^2}{4}\qquad\text{or}\qquad g=\frac{(q-1)q}{2}\, .
$$
Maximal curves with genus $(q-1)q/2$ or $(q-1)^2/4$ are unique (up to
$\fq$-isomorphism)(see \cite{rsti}, \cite[Theorem 3.1]{fgt}). In the latter
case, the proof of the uniqueness involves the study of the interplay
between the canonical divisor and the divisor $\cD:=|(q+1)P_0|$ on the
curve, $P_0$ being a $\fq$-rational point, via St\"ohr-Voloch's approach
to the Hasse-Weil bound (see \cite{sv}). Here we study
further the interplay between these divisors and we notice 
that the support of the $\fq$-Frobenius divisor associated to 
$\cD$ is contained in the union of the set of $\fq$-rational points and
the set of Weierstrass points of the curve (see Theorem \ref{t2.1}). 
This inclusion turned out to be an equality 
(Theorem \ref{t3.1}) in the case of maximal curves having classical
Weierstrass gaps. This answers Question \ref{q2.1} here at the
set-theoretical level and Examples \ref{e2.1} and \ref{e3.1} give
a positive answer to that question for a class of hyperelliptic maximal 
curves. We recall that the genus $g$ of a maximal curve over $\fq$
satisfies $g\ge q-n$, $n+1$ being the projective dimension of the linear
system $\cD$, and $g=q-n$ if the curve has classical Weierstrass gaps
\cite[Proposition 1.7(i)]{fgt}.  So while in \cite{fgt} we were interested
in maximal curves with high genus, here on the contrary we are mainly
interested on maximal curves with the smallest possible genus. 

We end up the paper by giving examples of maximal curves with classical
Weierstrass gaps.
\section{Maximal curves}\label{2} 
We use the following terminology and notations:
\begin{itemize}
\item A curve over $\fq$ is a projective geometrically irreducible
nonsingular algebraic curve defined over $\fq$. 
\item For $X$ a maximal curve over $\fq$ and $P_0\in X(\fq)$ we set
$$
\cD:= |(q+1)P_0|,\qquad n+1:= \dim(\cD)\, .
$$
\item For a $\fq$-linear system $\cL$ on $\xq$ we denote by $R^{\cL}$
(resp. $S^{\cL}=S^{(\cL,q^2)}$) the ramification
divisor (resp. the $\fq$-Frobenius divisor) associated to $\cL$;
the notation $j_i(P)$ 
(resp. $L_i(P)$) stands for the $i$th $(\cD,P)$-order
(resp. $i$th $\cD$-osculating space at $P$); see \cite{sv}.
\item We denote by $\cK=\cK_X$ the canonical linear system on $X$. 
Recall that $\cW_X:= \supp(R^{\cK})$ is the set of Weierstrass points of
$X$.
\item For $P\in X$, $m_i=m_i(P)$ denotes the $i$th non-gap at $P$, with  
$m_0(P):=0$, and $H(P)$ the Weierstrass semigroup at $P$. Recall that a
curve $X$ is classical iff $m_1(P)=g+1$ for each $P\not\in \cW_X$, $g$
being the genus of $X$.
\item $\frx$ denotes the Frobenius morphism on $X$ relative to $\fq$.
\end{itemize}
\begin{fund}[FGT, Corollary 1.2] For $\xq$ a maximal curve,  
$P_0\in X(\fq)$ and $P\in X$, we have the following linear equivalence:
$$
qP+\frx(P)\sim (q+1)P_0\, .
$$
\end{fund}
It follows that $n+1$ is independent of $P_0\in X(\fq)$, that
$m_{n+1}(P)=q+1$ for each $P\in X(\fq)$, and that
\begin{equation}\label{eq1}
m_0(P)=0<\ldots<m_n(P)\le q<m_{n+1}(P)\qquad \text{for each $P\in X$}\, .
\end{equation}
Therefore the following numbers are $(\cD,P)$-orders for $P\not\in X(\fq)$ 
\cite[Prop. 1.5(ii)]{fgt}:
\begin{equation}\label{eq2}
q-m_i(P)\qquad\text{for}\qquad i=0,1,\ldots,n\, .
\end{equation}
In addition, one can show that $m_n(P)=q$ for each $P\in X$  
\cite[\S2.3]{ft1}.
\begin{theorem}\label{t2.1} For $\xq$ a maximal curve, we have
$$
\supp(S^\cD)\subseteq\cW_X\cup X(\fq)\, .
$$
\end{theorem}
\begin{proof} Let $P\not\in \cW_X\cap X(\fq)$. Then $m_i(P)$ is
independent of $P$ and $\nu_i:=q-m_{n-i}(P)$ ($i=0,\ldots,n$) are the
$\fq$-Frobenius orders of $\cD$ \cite[\S2.2]{ft1}. Furthermore, by
(\ref{eq2}), there exists $I=I(P)\in \Z^+$ such that
\begin{equation*}
\nu_0<\ldots<\nu_{I-1}<j_I(P)<\nu_I<\ldots<\nu_n \tag{$*$}
\end{equation*}
are the $(\cD,P)$-orders. 
\begin{claim*}\quad $\frx(P)\in L_I(P)\setminus L_{I-1}(P)$.
\end{claim*}
\begin{proof} ({\it Claim}) For $i=0,1,\ldots,n$ let $u_{n-i}\in {\bar
\fq}(X)$, where $\bar\fq$ stands for the algebraic closure of $\fq$, be
such that $\df(u_{n-i})=D_i-m_i(P)P$, with 
$D_i\succeq 0$, $P\not\in\supp(D_i)$ and  
let $u\in {\bar \fq}(X)$ be such that $\df(u)=qP+\frx(P)-(q+1)P_0$ (cf.
the Fundamental Linear Equivalence). Then
\begin{equation*}
div(uu_{n-i})+(q+1)P_0=(q-m_i(P))P+\frx(P)+D_i\, .\tag{$**$}
\end{equation*}
By considering the morphism $\pi$ with homogeneous coordinates 
$v$ and $uu_{n-i}$, $i=0,1,\ldots, n$, where $v\in {\bar\fq}(X)$ is such
that 
$div(v)+(q+1)P_0=j_I(P)P+D_v$, with $D_v\succeq 0$, $P\not\in\supp(D_v)$,
we see
from \cite[Proof of
Thm. 1.1]{sv} and $(**)$ that $\frx(P)\in L_I(P)$. Now (loc. cit.)
if $\frx(P)\in L_{I-1}(P)$, then $\frx(P)\in \supp(D_v)$ and by the
Fundamental Linear Equivalence we
would have $q-j_I(P)=m_i(P)$ for some $i=0,1,\ldots,n$, a contradiction.
\end{proof}
To finish the proof of the theorem, notice that the claim and
$(*)$ imply the following linear relation
\begin{equation*}
D^{j_I(P)}\pi(P)=a\pi(\frx(P))+\sum_{i=0}^{I-1}a_iD^{\nu_i}\pi(P)\, 
,\tag{$***$}
\end{equation*}
where $a\neq 0, a_i\in {\bar \fq}$ and $D^j\pi(P)$ is the vector whose
coordinates are evaluations at $P$ of the Hasse derivatives (with
respect to a local parameter at $P$) of the homogeneous coordinates of the
morphism $\pi$ defined above. Now suppose that $P\in \supp(S^\cD)$. Then,
the following vectors would be linearly dependent
$$
\pi(\frx(P)), D^{\nu_0}\pi(P), D^{\nu_1}\pi(P),\ldots,D^{\nu_n}\pi(P)\, .
$$
From the linear relation in $(***)$ we then conclude that the following
vectors would be linearly dependent
$$
D^{j_I(P)}\pi(P), D^{\nu_0}\pi(P),
D^{\nu_1}\pi(P),\ldots,D^{\nu_n}\pi(P)
$$
and this contradicts the fact that the elements in $(*)$ are the
$(\cD,P)$-orders.
\end{proof}
\begin{remark}\label{r2.1} Recall that $X(\fq)\subseteq \supp(R^{\cD})$
\cite[Thm 1.4]{fgt} and that these sets may be different from each other 
\cite[Example 1.6]{fgt}. So in general we have that
$$
\supp(S^\cD)\subseteq \cW_X\cup X(\fq)\subseteq
\cW_X\cup \supp(R^\cD),
$$
where the last inclusion may be proper.
\end{remark}
\begin{remark}\label{r2.2} Suppose that $X$ is both maximal and classical. 
Then, by considering $P\not\in\cW_X$, from (\ref{eq1}) we have that 
$g=q-n$ and that the $\cD$-orders are $0,\ldots,n-1,\epsilon_n\ge n$ and  
$q$ (cf. \cite[Prop. 1.5(ii), Prop. 1.7]{fgt}). We
also have that the $\fq$-Frobenius orders of $\cD$
are $0,\ldots,n-1,q$ (\cite[\S2.2]{ft1}). Since \cite[p. 9]{sv}
$$
\deg(S^\cD)= \sum_{i=0}^{n}\nu_i (2g-2) + (q^2+n+1)(q+1)\, ,
$$
after some computations we find that
\begin{equation}\label{eq3}
\deg(S^{\cD})=(n+1)\#X(\fq)+\deg(R^{\cK_X})\, .
\end{equation}
This together with Theorem \ref{t2.1} suggest the following question 
\begin{question}\label{q2.1} For a classical maximal curve $\xq$, holds it
that
$$
S^{\cD}=(n+1)\sum_{P\in X(\fq)}\, P + R^{\cK_X}\ ?
$$
\end{question}  
\end{remark}
\begin{example}\label{e2.1} Here we are going to show that the
equality in Question \ref{q2.1} holds for certain hyperelliptic maximal
curves. Let $\xq$ be 
such a curve of genus $g>1$. By considering the unique linear system
$g^1_2$
on $X$ and the maximality of $X$ we see that $q\ge 2g$; furthermore, it
is well known that $X$ is classical. We set $\cW:= \cW_X$ and we restrict 
our attention to the case where one has 
$$
\cW\subseteq X(\fq)\qquad\text{and}\qquad \text{$q$ odd}\, .
$$
(The case $q$ even will we considered in Example \ref{e3.1}.) 
There are two types of $\fq$-rational points: either $P\in \cW$ or
$P\not\in \cW$.

Let $P\in \cW$. Setting $t:=q-2g$, the first
$(n+2)$-non-gaps of $X$ at $P$ are $0,2,\ldots,2g,2g+1,\ldots,2g+t,2g+t+1$
and hence,
by \cite[Prop. 1.5(iii)]{fgt}, the $(\cD,P)$-orders are
$0,1,\ldots,t+1,t+3,\ldots, 2g+t-1, 2g+t+1$. So \cite[Prop. 2.4(a)]{sv}
implies $v_P(S^{\cD})\ge q+(g-1)(g-2)/2$.

Let $P\in X(\fq)\setminus \cW$. As in the previous
case, here we find that $v_P(S^\cD)\ge n+1=q-g+1$.

Since $\#\cW=2g+2$ (here we use $q$ odd), then after some computations we
have that 
$$
\sum_{P\in X(\fq)}v_P(S^{\cD})\ge
\deg(S^\cD)=(2g-2)(\frac{(n-1)n}{2}+q)+(q^2+n+1)(q+1)\, ,
$$
hence that 
$$
S^\cD=\sum_{P\in \cW}\, (q+\frac{(g-1)(g-2)}{2})P+
\sum_{P\in X(\fq)\setminus\cW}\, (n+1)P\, .
$$
From this one concludes that the equality in Question \ref{q2.1} holds by
using the fact that the multiplicity of a Weierstrass point in the divisor
$R^\cK$ is $g(g-1)/2$.
\end{example}
\section{Certain maximal curves}\label{3}
The curves we have in mind in this section are maximal curves having
classical Weierstrass gaps, however we will consider a more general
setting. For a maximal curve $\xq$ let us consider the following conditions:
\begin{enumerate}
\item[(I)] For $Q_1, Q_2\not\in X(\fq)$,\ 
$H(Q_1)\cap [0,q]=H(Q_2)\cap [0,q] \Rightarrow H(Q_1)=H(Q_2)$.
\item[(II)] For $Q\not\in \cW_X$, $m_i(Q)=q-n+i$, $i=1,\ldots, n$.
\end{enumerate}
By (\ref{eq1}) each maximal curve with $g=q-n$ (e.g. 
a classical maximal curve) satisfies Condition (I). Other examples are  
provided by maximal curves $\xq$ with $\cW_X\subseteq X(\fq)$; in
this case $\cW_X=X(\fq)$ whenever $g>q-n$ \cite[Corollary 2.3]{ft1}. 

Condition (II) is satisfied by classical maximal curves, by the
Hermitian curve and by some curves covered by this curve (see \cite{gv}). 

For $\xq$ a maximal curve, denote by $\tilde m_i$ the $i$th non-gap
at $P\not\in\cW_X$. If $X$ satisfies Condition (II), then, by  
(\ref{eq2}) and \cite[Thm. 1.4(i)]{fgt}, the $\cD$-orders are\\
$0,\ldots,
n-1, \epsilon_n\ge n$ and  
$\epsilon_{n+1}=q$; furthermore, by \cite[\S2.2]{ft1}, the $\fq$-Frobenius
orders are $0,\ldots,n-1$ and $q$.
\begin{proposition}\label{p3.1} For a maximal curve $\xq$ satisfying
Condition (I), one has
$$
\cW_X\setminus \supp(R^\cD) \subseteq \supp(S^\cD)\, .
$$
\end{proposition}
\begin{proof} We first notice that, for $Q\not\in X(\fq)$, Condition (I) 
implies
$$
Q\in \cW_X\quad \Leftrightarrow\quad \{m_1(Q),\ldots,m_n(Q)\}\neq\{\tilde
m_1,\ldots,\tilde m_n\}\, ;
$$
we also notice that $m_i(Q)\le \tilde m_i$ for each $i$. 
Now let $Q\in \cW_X\setminus \supp(R^\cD)$ and let $k\in [1,n-1]$ be such
that $m_i(Q)=\tilde m_i$ for $1\le i<k$ and $m_k(Q)<\tilde m_k$. Then, by
(\ref{eq2}), the $\cD$-orders (which are also the $(\cD,Q)$-orders since
$Q\not\in\supp(R^\cD)$) are
\begin{align*}
\epsilon_0= & 0=q-\tilde m_n< \epsilon_1=1=q-\tilde
m_{n-1}\ldots<\epsilon_{n-k}=q-\tilde m_k <\\
  & \epsilon_{n-k+1}=q-m_k(Q)<\epsilon_{n-k+2}=q-\tilde
m_{k-1}<\ldots<\epsilon_{n+1}=q=q-\tilde m_0
\end{align*}
so that $m_i(Q)=\tilde m_i$ for $k+1\le i\le n$. As in the proof of the
claim in Theorem \ref{t2.1}, we conclude that (for $J=n-k$)
$$
\frx(Q)\in L_J(Q)\setminus L_{J-1}(Q)\, .
$$
We note also that $\epsilon_{J+1}$ is the $\cD$-order one should take out
to get the $\fq$-Frobenius orders of $\cD$, hence $\nu_i=\epsilon_i$ for
$i\le J$. We then conclude that the vectors
$$
\pi(\frx(Q)), D^{\nu_0}\pi(Q),\ldots,D^{\nu_J}\pi(Q)
$$
are linearly dependent and this finishes the proof.
\end{proof}
\begin{lemma}\label{l3.1} Let $\xq$ be a maximal curve satisfying
both Conditions (I) and (II) and let $P\in \cW_X\setminus X(\fq)$ with
$j_{n-1}(P)=n-1$. Then $m_1(P)=q-j_n(P)$.
\end{lemma}
\begin{proof} By (\ref{eq2}) we have that $q-m_1(P)$ is a $(\cD,P)$-order 
with $q-m_1(P)\le j_n(P)$. If $q-m_1(P)<j_n(P)$ we would have $q-m_1(P)\le
n-1$ so that $m_i(P)=q-n+i$ for $i=1,\ldots, n$ (see (\ref{eq1})). 
Consequently by Condition (II), $H(P)\cap [0,q] = H(Q)\cap [0,q]$, for 
$Q\not\in
\cW_X$ and hence $H(P)=H(Q)$ by Condition (I), i.e. $P$ is not a
Weierstrass point, a contradiction.
\end{proof}
\begin{remark}\label{r3.1} Suppose that $X$ satisfies both Conditions 
(I) and (II). Let \\
$P\in \cW_X\setminus \supp(R^\cD)$. Then Lemma
\ref{l3.1} implies $m_1(P)=q-\epsilon_n$ and from the proof of Proposition
\ref{p3.1}, we have $m_i(P)=q-n+i$, $i=2,\ldots, n$. Then 
$$
n\le \epsilon_n\le \frac{q+n-2}{2}\, ,
$$
where the last inequality follows from the fact that $2m_1(P)\ge m_2(P)$. 
Next we state a criterion to ensure that $\epsilon_n=n$, namely
$$
\cW_X\setminus \supp(R^\cD)\neq\emptyset\quad \text{and}\quad p:= {\rm
char}(\fq)\ge
g\qquad \Rightarrow\qquad \epsilon_n=n\, .
$$
Indeed if $\epsilon_n>n$, by the $p$-adic criterion \cite[Corollary 
1.9]{sv}, we
would have $\epsilon_n\le q-p$ so that $\epsilon_n\le n+g-p$ (see
(\ref{eq1})), i.e. $\epsilon_n=n$, a contradiction.
\end{remark}
\begin{lemma}\label{l3.2} Let $\xq$ be a maximal curve satisfying
both Conditions (I) and (II) and let $P\in X\setminus X(\fq)$ with
$j_{n-1}(P)=n-1$. Then the following statements are equivalent:
\begin{enumerate} 
\item $q-j_n(P)\in H(P)$.
\item $P\in \supp(S^\cD)$.
\end{enumerate}
\end{lemma}
\begin{proof} $(1)\Rightarrow (2):$ If $q-j_n(P)\in H(P)$, then as in the
proof of the claim in Theorem \ref{t2.1} we have 
$\frx(P)\in L_{n-1}(P)$ and hence it belongs to  
$\supp(S^\cD)$ since $0,\ldots, n-1$ are the $\fq$-Frobenius
orders of $\cD$ and $j_{n-1}(P)=n-1$.
\smallskip

\noindent $(2)\Rightarrow (1)$ By Theorem \ref{t2.1} we have that $P\in
\cW_X$ and the result follows from Lemma \ref{l3.1}.
\end{proof}
\begin{corollary}\label{c3.1} Let $\xq$ be a maximal curve
satisfying both Conditions (I) and (II). For 
$P\in X$ the following statements are equivalent:
\begin{enumerate}
\item $P\not\in \supp(S^\cD)$.
\item The $(\cD,P)$-orders are $0, 1,\ldots,n-1, j_n$ and $j_{n+1}=q$ 
with $q-j_n\not\in H(P)$.
\end{enumerate}
\end{corollary}
\begin{proof} By Lemma \ref{l3.2} we just need to show  that 
$(1)\Rightarrow (2)$. Since $0,\ldots,n-1$ are $\fq$-Frobenius
orders and in view of \cite[Thm. 1.4(ii)]{fgt}, we see
that $j_{n-1}(P)>n-1$ or 
$j_{n+1}(P)=q+1$ imply that $P\in\supp(S^\cD)$. So the $(\cD,P)$-orders
are as stated in (2) and the result follows again from Lemma \ref{l3.2}
\end{proof}
From Theorem \ref{t2.1}, Lemma \ref{l3.1} and Lemma \ref{l3.2}, we obtain
\begin{corollary}\label{c3.2}
Let $\xq$ be a maximal curve satisfying both Conditions (I) and (II)
and $P\in X\setminus X(\fq)$ with $j_{n-1}(P)=n-1$. The following
statements are equivalent:
\begin{enumerate}
\item $P\in \supp(S^\cD)$.
\item $P\in \cW_X$.
\item $m_1(P)=q-j_n(P)$.
\item $q-j_n(P)\in H(P)$.
\end{enumerate}
\end{corollary} 
Now we can state the main result of this section:
\begin{theorem}\label{t3.1} For a maximal curve satisfying 
both Conditions (I) and (II), we have 
$$
\supp(S^\cD)=\cW_X\cup X(\fq)\, .
$$
\end{theorem}
\begin{proof} From Theorem \ref{t2.1} and the fact that $X(\fq)\subseteq
\supp(S^\cD)$ it is enough to show that
$$
P\in \cW_X\setminus X(\fq)\quad \Rightarrow\quad P\in \supp(S^\cD)\, .
$$
If $j_{n-1}(P)>n-1$, then $P\in \supp(S^\cD)$ (see the proof of Corollary
\ref{c3.1}). So let now $j_{n-1}(P)=n-1$. Then again $P\in \supp(S^\cD)$  
as follows from Corollary \ref{c3.2}.
\end{proof}
\begin{example}\label{e3.1} We complement Example \ref{e2.1} by 
considering hyperelliptic maximal curves of genus bigger than 1 over 
$\fq$ with $q$ even. We are going to show that these curves also satisfy 
the equality in Question \ref{q2.1}. So let $X$ be a such curve. By
\cite[Proposition 1.7(ii)]{fgt}, $\#\cW_X=1$; say $\cW_X=\{Q\}$. Then for 
at least $(\#X(\fq)-1)$ $\fq$-rational points $P\in X$ we have  
$v_P(S^\cD)=n+1$, as follows
from the computations in Example \ref{e2.1} and \cite[Proposition
2.4(a)]{sv}. Furthermore $\supp(S^\cD)=\{Q\}\cup X(\fq)$ by Theorem
\ref{t3.1}. Next we compute $v_Q(S^\cD)$ by using Eq. (\ref{eq3}). We
consider two cases according $Q$ is $\fq$-rational or not.

If $Q\in X(\fq)$, from (\ref{eq3}) we have 
$$
v_Q(S^\cD)=\deg(S^\cD)-(n+1)(\#X(\fq)-1)=\deg(R^\cW_X)+n+1\, .
$$

If $Q\not\in X(\fq)$, from (\ref{eq3}) we have
$$
v_Q(S^\cD)=\deg(S^\cD)-(n+1)\#X(\fq)=\deg(R^\cW_X)\, .
$$
From these computations follow the equality in Question \ref{q2.1} for hyperelliptic maximal curves 
over $\fq$ with $q$ even.
\end{example}
\section{Examples}\label{4}
From \cite{rsti}, the unique maximal curve over $\fq$ of genus $q(q-1)/2$
is
the Hermitian curve in $\P^2(\bar\fq)$ defined by
$$
Y(Y^{q-1}+1)=X^{q+1}\, .
$$
Let $\pi:\P^2(\bar\fq)\to \P^2(\bar\fq)$ be the morphism over $\fq$ given
by\\ 
$(x:y:1)\to (y^{q-1}:x^{(q^2-1)/m}:1)$ with $m$ a divisor of $(q^2-1)$. 
Then the nonsingular model of $\pi(X)$ is a maximal curve over $\fq$
\cite[Proposition 6]{lachaud} and $\pi(X)$ 
is defined by 
\begin{equation}\label{eq4}
W^m=Z(Z+1)^{q-1}\, .
\end{equation}
By the Riemann-Hurwitz relation, the genus $g$ of this curve satisfies
$$
g=\frac{m-\delta}{2},\qquad \text{where $\delta=\gcd{(m,q-1)}$}\, .
$$
These examples are the ones in \cite[Corollary 4.9]{gstix} (see also
Lang's \cite[Ch. I,\S7]{lang} and the references therein). Suppose one is
interested in genus 4 maximal curves of the type above. So $m-\delta=8$
and since $\delta$ divides $m$, we have that $\delta=1,2,4$ or 8. 
As an example, let $\delta=1$ and hence $m=9$. Since $m=9$ divides 
$(q^2-1)$ and moreover $\delta=1=\gcd{(m,q-1)}$, we must have that
$m$ divides $(q+1)$. So the prime power $q$ must be chosen in the
following congruence
class:  
$$
q\equiv -1\pmod{9}\, .  
$$
With the above reasoning one obtains the following table which gives for a
fixed genus $g$ ($1\le g\le 7$) maximal curves over $\fq$ arising from
curves of type (\ref{eq4}). 
%
%
\bigskip

\begin{center}
\begin{tabular}{|c|c|c||c|c|c|}
\hline
{\rm Genus} & $m$ & $q$ & {\rm Genus} & $m$ & $q$ \\
\hline
\hline
$g=1$ & 3 & $q\equiv -1 \pmod 3$ & $g=5$ & 11 & $q\equiv -1\pmod{11}$\\
\hline
$g=1$ & 4 & $q\equiv -1 \pmod 4$ & $g=5$  & 12 & $q\equiv -1 \pmod{12}$\\
\hline
$g=2$ & 5 & $q\equiv -1 \pmod 5$ & $g=5$ & 15 & $q\equiv -4\pmod{15}$\\
\hline
$g=2$ & 6 & $q\equiv -1 \pmod 6$ & $g=5$ & 20 & $q\equiv 11\pmod{20}$\\
\hline
$g=2$ & 8 & $q\equiv 5 \pmod 8$ & $g=6$ & 13 & $q\equiv -1\pmod{13}$\\
\hline
$g=3$ & 7 & $q\equiv -1 \pmod 7$ & $g=6$ & 14 & $q\equiv -1\pmod{14}$\\
\hline
$g=3$ & 8  & $q\equiv -1 \pmod 4$ & $g=6$ & 15 & $q\equiv 4\pmod{15}$\\
\hline
$g=3$ & 12 & $q\equiv 7 \pmod{12}$ & $g=6$ & 24 & $q\equiv 13\pmod{24}$\\
\hline
$g=4$ & 9 & $q\equiv -1 \pmod 9$ & $g=7$ & 15 & $q\equiv -1\pmod{15}$\\
\hline
$g=4$ & 10 & $q\equiv -1 \pmod{10}$ & $g=7$ & 16 & $q\equiv -1\pmod{8}$\\
\hline
$g=4$ & 12 & $q\equiv 5 \pmod{12}$ & $g=7$ & 21 & $q\equiv 8 \pmod{21}$\\
\hline
$g=4$ & 16 & $q\equiv 9 \pmod{16}$ & $g=7$ & 28 & $q\equiv 15\pmod{28}$\\
\hline
\end{tabular}
\end{center}
\bigskip

%
By the Dirichlet theorem there are infinitely many prime numbers in each
of the congruence classes above. Notice that the classicality of the curve
is assured as soon as ${\rm char}(\fq)>2g-2$.

We finish this section with a more restricted class of classical maximal
curves . We consider maximal curves over $\fq$ where each rational point
is not a Weierstrass point, so that $v_P(R^\cD)=1$ for each $P\in X(\fq)$,  
and where
$$
X(\fq)=\supp(R^\cD)\, .
$$
Then from 
$$
\deg(R^\cD)=(\frac{n(n+1)}{2}+q)(2g-2)+(n+2)(q+1)=\#X(\fq)=(q+1)^2+q(2g-2)
$$
we have that $n(n+1)(g-1)=(q+1)(g-1)$ and hence
$$
g=1\qquad\text{or}\qquad q=n^2+n-1\, .
$$
The only example we know of a curve with $g>1$ in this restricted class
of maximal curves is the one over $\F_{25}$ of genus $3$ listed
by Serre in \cite[Section 4]{serre}; in this example $n=2$ (cf.
\cite[Example 2.4(i)]{fgt}). 

\end{document}